\documentstyle[pb-diagram,lamsarrow,pb-lams,amssymb,amsfonts,amsthm]{amsart}

\def\mystretch{1.2}
\def\baselinestretch{\mystretch}

\def\baselinestretch{\mystretch}

\newtheorem{theorem}{Theorem}[section]

\newtheorem{proposition}[theorem]{Proposition}
\newtheorem{lemma}[theorem]{Lemma}

\theoremstyle{definition}
\newtheorem{definition}[theorem]{Definition}

\theoremstyle{remark}
\newtheorem{remark}[theorem]{Remark}

\def\abstract{\noindent \vspace{0.3cm}\large{\bf Abstract}\\ 
\small\def\baselinestretch{1}\normalsize}

\newcommand{\W}{{\mathcal {W}}} 
 
\newcommand{\A}{{\mathbb {A}}}
\newcommand{\C}{{\mathbb {C}}}

\newcommand{\R}{\mathbb {R}}

\newcommand{\dispt}{\displaystyle\frac{i}{t}}

\begin{document}

\title{How to calculate the Fedosov star--product\\
{\tt (Exercices de style)}
}
\author {Olga Kravchenko}
\address{Institut Girard Desargues UMR 5028, Université Lyon-I, 
     43, boulevard du 11 novembre 1918,
        69622 Villeurbanne Cedex, France, ok@@alum.mit.edu}

\thanks{Presented in RENCONTRES MATHÉMATIQUES DE
                    GLANON (Bourgogne, France)  1998\\
\indent http://www.u-bourgogne.fr/glanon/\\}

\maketitle

\begin{flushright} {\it To the memory of Moshé Flato}
\end{flushright} 

\vspace{.5cm}

\begin{abstract}

\noindent  This is an expository note on Fedosov's construction of deformation
quantization. Given a symplectic manifold and a connection on it, we show how
to calculate the star-product step by step.

\noindent  We draw simple diagrams to solve the recursive 
equations for the Fedosov connection and for flat sections
of the Weyl algebra bundle corresponding to functions.

\noindent   We also reflect on the differences of symplectic and Riemannian geometries.

\end{abstract}

\vspace{1cm}

\noindent {\bf Key--words:}  deformation quantization, 
 Weyl algebra, Fedosov quantization\\
\noindent {\bf  AMS classification (2000):} 53D55

%%%%%%%%%%%%%%%%%%%%%%%%%%%%%%%%%%%%%%%%%%%%%%%%%%%%%%%%%

\section{Deformation quantization of a symplectic manifold}

\subsection{Fedosov's idea: Koszul--type resolution.\\}
We consider a deformation quantization 
 of a symplectic manifold $(M, \omega_0)$ as 
 a deformation of an 
algebra of smooth functions on $M$ in the direction of the Poisson
 bracket  \cite{Flato}. 

\begin{definition}
Deformation quantization of a symplectic manifold $(M, \omega_0)$ is 
an associative algebra structure on ${\A}= C^\infty (M) [[t]]$ over $\C[[t]],$ 
called  a $\ast$-product,
such that for any $a = a(x,t) = \sum_{k=0}^\infty t^k a_k(x) \ \mbox{and} \ 
  b = b(x,t) = \sum_{k=0}^\infty t^k b_k(x), \ a_k(x),  b_k(x) 
\in C^\infty (M) $ 
\begin{enumerate}
\item The product $\ast$ is local, that is in the $\ast$-product 
$
a(x, t) \ast b(x,t) = \sum_{k=0}^\infty t^k c_k(x),
$
the coefficients $ c_k(x)$ depend only on $a_i, b_j$ and their derivatives 
$\partial^\alpha a_i, \partial^\beta b  $ with 
$i+j+ |\alpha| + |\beta| \leq k.$
\item It is a formal deformation of the commutative 
algebra $C^\infty(M): \ c_0 (x) = a_0 (x) b_0(x).$ 
\item Let $\{ \cdot, \cdot \}$ be  the Poisson 
bracket of functions, given by a bivector field  
dual  to the form $\omega_0.$ 
There is a correspondence principle:
\[
[a,b] : = \dispt ( a \ast b- 
 b \ast  a) =\{a_0(x), b_0(x)\} + t \ r(a,b),
\]
where $r(a,b) \in \A.$
\item There is a unit:
$a(x, t) \ast 1 = 1 \ast a(x, t) = a(x, t).$
\end{enumerate}
\end{definition}

Fedosov  found a geometric way to perform the deformation 
quantization \cite{F,F1} (also see \cite{Le} for a comprehensive 
exposition).  
The following idea lies behind Fedosov's construction  --- a
Koszul--type resolution is considered for $C^{\infty}(M)[[t]]:$
 \begin{equation}
\label{complex}
 C^{\infty}(M)[[t]] \stackrel{ Q} \to A^0 
\stackrel{D} \to A^1
\stackrel{D} \to A^2 \stackrel{D} \to \ldots. 
 \end{equation}
It means that cohomology groups  
of this complex are all zero, and in particular  
we get  $ C^{\infty}(M)[[t]]: = Im Q \cong Ker D.$  
If each term of the resolution has a noncommutative associative structure 
and moreover $D$ respects this structure, then it provides the space 
$ C^{\infty}(M)[[t]] $ with a new associative 
noncommutative product. 
\noindent  Namely,
let $\circ : A^0 \otimes A^0  \to A^0$ be such a product on $A^0.$ 
Then we get a product on $ C^{\infty}(M)[[t]] $ as follows:
\begin{equation}
\label{eq:ast}
a \ast b = Q^{-1} \bigl( Q (a) \circ Q (b) \bigr). 
\end{equation}
($D$ respects the product, so since 
$D$ is zero on $Q (a)$ and  $Q (b)$, it must be  zero on  $Q (a) \circ Q (b),$
so the product is in the kernel of $D$, that is in the  image of $Q$, and its preimage 
 $Q^{-1}$ is well defined.)
The product on $A^0$
should also verify certain properties to 
certify the axioms of the deformation quantization, so 
it is a very special resolution.

 Fedosov constructs  such a resolution by using the differential 
forms on the manifold with values in the  
Weyl algebra bundle (with Moyal-Vey fibrewise product) \cite{F}. 
The main step then is to find  the  differential on it, $D,$  which 
would respect the algebra structure. This
differential is called Fedosov connection and is obtained by an iteration
procedure from a torsion free 
symplectic connection on the  manifold. 

\subsection{Weyl algebra of a vector
space.\\}
 
Let $E$ be a vector space with a non--degenerate skew-symmetric form $\omega$.
  The algebra of polynomials on $E$ is the algebra of symmetric
powers of $E^*$, $S(E^*)$, and it has a skew-symmetric form on it 
which is dual to $\omega$.
Let $e$ be a point in $E$ and $\{e^k\}$ denote its linear coordinates in
$E$ with respect to some fixed basis.  Then $\{e^k\}$ define a basis in $E^*$. 
Let $\omega^{kl}$ be the
matrix for the skew-symmetric form on $E^*$. Let us consider the power
series in $t$ with values in $S(E^*)$: 

\begin{definition}
The Weyl algebra $W(E^*)$ of a vector space $E^*$ 
is an associative algebra
 \[ 
W(E^*) = S(E^*) [[T]]: \ \ \ 
 a(e,t) = \sum_{k \geq 0} a_k(e) t^k,
\] 
 with  the product structure given by the Moyal--Vey product:
	\begin{equation}
	\left.
 a\circ b (e,h)= \exp \  \{- \frac{it}{2} \omega ^{kl} \frac{\partial  }
{\partial
 x^k} \frac{\partial  }{\partial  z^l}\} \ \  a(x, t) \ 
b(z, t)\  \right| _{x=z=e}.
							 \label{eq:Moyal}
	\end{equation}
 \end{definition}

\noindent 
The Lie bracket is defined with respect to this product.
We can look at this algebra as at a completion of 
the universal enveloping algebra of the 
Heisenberg algebra on $E^* \oplus t \C $, namely, the algebra 
with relations
\begin{equation}
\label{eq:homo}
e^k \circ e^l - e^l \circ e^k= - i t \omega^{kl}
\end{equation}
where $\omega ^{kl} = \omega (e^k, e^l)$ defines  a  Poisson bracket
on $E^*$.

\noindent Let us consider the product of the Weyl algebra and the 
exterior algebra of the space $E^*$: $W(E^*)\otimes \Lambda E^*$.
Let $dx^k$ be the basis in $\Lambda E^*$ corresponding to $e^k$ in $W(E^*)$.
 
\noindent There is a decreasing  filtration on the Weyl 
algebra $W(E^*)$: $
W_0 \supset W_1 \supset W_2 \supset ...$  given by the degree of
generators. The generators 
$e$'s have degree 1 and $t$ has degree 2: 
\[
W_p = \{ \text{elements with degree}  \geq p \}. 
\]

\noindent One can define a grading on
$W$ as  follows
\[
gr_i W = \{ \text{elements with degree}  = i \}. 
\]
it is isomorphic to 
$W_i/W_{i+1}$.
One can see  that  the product (\ref{eq:Moyal}) preserves the grading
(since  the relation (\ref{eq:homo}) is homogeneous).

\begin{definition}
\label{def:gr}
An operator on  $W(E^*)\otimes \Lambda E^*$ 
is said to be of degree $k$ if it maps 
$W_i\otimes\Lambda E^*$ to $W_{i+k}\otimes\Lambda E^*$ for all $i$.
\end{definition}

\noindent 
Such an operator defines maps $ gr_i W\otimes\Lambda E^*$ to
$gr_{i+k}W\otimes\Lambda E^*$ for all $i$.

\begin{definition} 
Derivation on  $W(E^*)\otimes \Lambda E^*$ is 
a linear operator which satisfies the Leibnitz rule:
\[
D (ab) = (Da) b + (-1)^{\tilde{ a}\tilde {D}}a (Db)
\]
\end{definition}
\noindent 
where $\tilde{a }$ and $\tilde{D }$ are  corresponding degrees. 
It turns out that all $\C[[t]]$--linear derivations are inner \cite{D}.
\begin{lemma}
 \label{l:inner}
 Any linear derivation $D$ 
on $W (E^*) \otimes \Lambda E^*$ is inner, 
 namely there exists such $v\in W(E^*)$ so that 
$D a = \dispt [v,a]$ for any $a \in W (E^*)$
\end{lemma}
\begin{pf} 
Indeed, for  the generators
$\frac{\partial}{\partial e^k} a 
= \frac{i}{2 t}[ \omega_{kl} e^l, a].$  
So for any derivation one can get a formula:
$D a = \dispt [ \frac{1}{2} \omega_{kl} e^k D e^l , a].$  
\end{pf}
One can define two natural  operators on the algebra 
$W(E^*)\otimes \Lambda E^*$: $\delta$ and 
$\delta^*$ of degree $-1$ and $1$ correspondingly, such that 
$\delta$ is the lift of the ``identity'' operator 
\[
u: e^k \otimes 1 \to 1 \otimes dx^k
\]
and  $\delta^*$ is the lift of its inverse.
On monomials  $e^{i_1}\otimes...\otimes e^{i_m}\otimes dx^{j_1} \wedge ...
\wedge dx^{j_n} \in W^m(E^*)\otimes \Lambda^n E^*$ 
$\delta$ and 
$\delta^*$ can be written as follows:
\begin{eqnarray*}
\lefteqn{
\delta: \ e^{i_1}\otimes...\otimes e^{i_m}
\otimes dx^{j_1} \wedge ...\wedge dx^{j_n}
\mapsto}\\
& &\sum _{k=1}^m
e^{i_1}\otimes...\widehat{e^{i_k}}...\otimes
e^{i_m}\otimes dx^{i_k} \wedge dx^{j_1} \wedge ...\wedge dx^{j_n}\\
\lefteqn
{\delta^*: \  e^{i_1}\otimes...\otimes e^{i_m}
\otimes dx^{j_1} \wedge ...dx^{j_n}
\mapsto}\\
& &\sum _{l=1}^n (-1)^l e^{j_l} \otimes
e^{i_1}\otimes...\otimes e^{i_m}\otimes dx^{j_1} \wedge ...
\widehat{dx^{j_l}} ... \wedge dx^{j_n}.
\end{eqnarray*}
\begin{lemma}
\label{l:delta}
Operators $\delta$ and $\delta^*$ have the following properties:
\[
\delta a = dx^l  \frac{\partial a }{\partial  e^l} =
\ [ -\frac{i}{t} \omega_{kl} \ e^k dx^l, \ a], \quad
\delta ^* a = y^l \iota_{\frac{\partial }{\partial  x^l}} a, \quad
\delta^2 = {\delta^*}^2 = 0
\]
On monomials $e^{i_1}\otimes...\otimes e^{i_m}\otimes dx^{j_1} \wedge ...
\wedge dx^{j_n} $ from $ gr_m W(E^*)\otimes \Lambda^n E^*$ 
\[
\delta \delta^* + \delta^* \delta = (m+n) Id,
\] 
where $Id$ is the identity 
operator.
Any element $a \in gr_m W(E^*)\otimes \Lambda^n E^*$ has a decomposition:
\[
a = \frac{1}{m+n}(\delta \delta^* a + \delta^* \delta a)  + a_0.
\]
where   $a_0$ is  a projection  of $a \in W(E^*)\otimes \Lambda E^*$ to 
the center of the algebra, that is 
the summands in  $a$ which  do not contain $e$-s.
\end{lemma}

%%%%%%%%%%%%%%%%%%%%%%%%%%%%%%%%

\subsection{Symplectic connections (symplectic covariant derivatives).\\}

\label{sec:sharp}

The term  {\em symplectic connection} in this section is in fact a
 {\em symplectic covariant derivative}.

Let us consider connections on a manifold $M$.
	\begin{proposition}
Let $\omega$ be a skew-symmetric 2-form on ${\mathcal T}M$.
Let  $\nabla$ be a torsion free connection preserving this form.
Then $\omega$ is necessarily closed.
	\end{proposition}
\begin{pf} The skew-symmetry of $\omega$ is the following condition: 
$\omega (X,Y) = - \omega (Y,X)$.  
The connection $\nabla$ is torsion--free when 
$\nabla_X Y -\nabla_Y X  = [X,Y]$. 
Suppose such $\nabla$ exists. The connection $\nabla$
 preserves the form
$\omega$ when $\nabla \omega = 0$. 
This means that for all $X,Y, Z \in {\mathcal T}M$:
\begin{equation}
\nabla_X (\omega ( Y, Z)) = \omega (\nabla_X Y , Z) +\omega ( Y ,\nabla _X Z) 
				\label{eq:nabla}
\end{equation}
Since $ \omega ( Y, Z)$ is a function 
 $\nabla_X ( \omega ( Y, Z)) = X \omega  ( Y, Z)$.
Then, 
\begin{eqnarray*}
X \omega  ( Y, Z) & - & Y \omega  ( X, Z) \ \ + \ \  Z \omega ( X, Y) \\
 & = & \omega (\nabla _X Y , Z) - \omega (\nabla _X Z, Y) -
 \omega (\nabla _Y X, Z) \\
 &  &\mbox{}+ \omega (\nabla _ Y Z, X) +
 \omega (\nabla _Z X, Y) - \omega (\nabla _Z Y, X)\\
& = & \omega ([X, Y], Z) - \omega ([X,Z], Y) + \omega ([Y, Z], X)
\end{eqnarray*}
which is exactly the condition $d \omega = 0$.
\end{pf}

\begin{remark} Here we want to make an analogy with a Riemannian case. 
The Riemannian metric is a symmetric two-form and there is
a unique torsion free connection compatible with it (the
Levi--Civita connection). 

\noindent In symplectic case we deal with a skew-symmetric form.
There are many different torsion free connections preserving the form if 
the form is closed.

\noindent The statement of uniqueness of Levi--Civita
connection in the Riemannian case is substituted by the 
requirement for the form to be closed in the skew-symmetric setting:
\begin{itemize}
\item 
{\large\sc Symmetric:} Torsion--free compatible connection 
always exists and unique  (Levi--Civita connection).
\item 
{\large\sc  Skew-symmetric:} Torsion--free compatible connection 
exists if the form is closed and, in general, not unique.
\end{itemize} 
\end{remark}

Here we are mostly interested in the case when $M$ is a 
symplectic manifold, that is  when there is a symplectic form
$\omega$ on $M$
(a closed and nondegenerate $2$-form  on ${\mathcal T}M$).
\begin{definition}
\label{def:con}
A connection which preserves  a symplectic form is called a 
symplectic connection.
\end{definition}
Any connection on a symplectic manifold gives rise 
 to a symplectic connection:

\begin{proposition} \cite{L} \cite{MRR}.
Let $(M, \omega) $ be a symplectic manifold.
Then for every 
connection $\nabla$ there exists a three--tensor $S$, such that 
\[
\tilde{\nabla} = \nabla + S
\]
is a symplectic connection.\\ 
Then for $X,Y \in {\mathcal T}M$
\[
\hat{\nabla}_X Y = \tilde{\nabla}_X Y - \frac{1}{2} Tor (X,Y)
\]
defines a torsion--free connection compatible with the form $\omega$. Here 
$2$--form $Tor$ is the  torsion of $\tilde{\nabla}$
\[
Tor (X,Y) = \tilde{\nabla} _X Y -   \tilde{\nabla} _Y X -
  \tilde{\nabla} _{[X,Y]}
\]
\end{proposition}
Then $S$ is defined as follows:
\[
S_X Y = \frac{1}{2} \{ (\nabla_X \omega )( Y, .) \}^{\sharp},
\]
where $\sharp:  {\mathcal T}^*M \to {\mathcal T}M $ is  the inverse  
to $\flat$, given by :
\begin{eqnarray*}
\flat:{\mathcal T}M \to {\mathcal T}^*M\\
u^\flat = \omega (u, .) \ for \ \ u \in {\mathcal T}M.
\end{eqnarray*}
% Let ${\cal A}^* (M)$ denote the algebra of differential forms on $M$.
Symplectic connections form an affine space with the  associated vector 
space    ${\mathcal A}^1(M,sp(2n)),$ the 
Lie algebra  $sp(2n)$--valued one-forms on $M$.

\subsection{Weyl algebra bundle and Fedosov's theorem.\\}
%%%%%%%%%%%%%%%%%%%%%%%%%%%%%%%%%%%
Let  $M^{2n}$ be a symplectic manifold with a symplectic form 
$\omega$. In local 
coordinates at a  point $x$:
      \[
\omega = \omega_{kl} d x^k \wedge dx^l.
\]
The symplectic form on a manifold $M$ defines a Poisson  
bracket on functions on $M$.
For any two functions $a, b \in C^{\infty}(M)$:
	\begin{equation}
	\{a,b \} = \omega^{kl} \frac{\partial  a}{\partial  x^k} 
\frac{\partial  b}{\partial  x^l}    \label{eq:poisson}
	\end{equation}
where $(\omega^{kl})=(\omega_{kl})^{-1}$ (as matrix coefficients).

We can define the bundle of Weyl algebras  ${\cal W} \to M$, with 
the fibre at a point $x \in M$
being the Weyl algebra of the vector space ${\mathcal T}^*_x M$. Let  
$\{e^1, ... e^{2n}\}$ be $2n$ 
generators in  ${\mathcal T}^*_x M$, 
corresponding to $dx^k$. The form $\omega^{kl}$ defines a 
pointwise Moyal--Vey product. 

\noindent

The filtration and the grading in ${\cal W}$ are inherited 
from $W({\mathcal T}^*_xM)$ at each point $x \in M$. We denote 
by  ${\cal W}^k$ the $k$-th graded component in ${\cal W}$: 
 \[
{\cal W} = \oplus_i {\cal W}^k
\]
 A symplectic connection, $\nabla$, satisfying (\ref{eq:nabla})
can be naturally lifted to act on any symmetric power 
of the cotangent bundle (by the Leibniz rule). Moreover, since 
the cotangent bundle ${\mathcal T}^*M \cong {\cal W}^1,$ we 
can lift $\nabla$ to be an operator on sections 
$\Gamma (M,{\cal W}^k)$ with values in 
$\Gamma (M, {\mathcal T}^*M \otimes {\cal W}^k )$.
By abuse of notations this 
operator is also called $\nabla$. 

\noindent The connection $\nabla$ preserves the grading, in other words 
it is an operator of degree zero.
It is clear that in general this connection is not flat: $\nabla ^2 \not= 0$.
Fedosov's idea is that for ${\W}_M$ bundle one can add to 
the initial symplectic connection some 
operators not preserving the  grading
so that the sum gives a flat connection on the Weyl bundle.

\begin{theorem}(Fedosov.)
\label{Fedosov}
There is a unique set of operators 
$r_k : \Gamma (M,{\cal W}^k) \to
\Gamma (M, {\mathcal T}^*M \otimes {\cal W}^{i+k})$ such that 
\begin{equation}
	D = - \delta + \nabla + r_1 + r_2 + \ldots 
				  \label{eq:connection}
	\end{equation}
\[
D^2 = 0, \ \mbox{and} \ \delta^* r_i =0.
\]
There is a one-to-one correspondence between formal series in $t$ 
with coefficients in smooth functions $C^\infty(M)$
and  horizontal sections of this connection:
\begin{equation}
\label{eq:Q}
Q: C^\infty(M)[[t]] \to  \Gamma_{\text{flat}} (M,{\cal W}_M).
\end{equation}
\end{theorem}
\noindent Main idea of the proof is to use the following complex:
 \begin{equation}
\label{comp:del}
0 \stackrel{ } \to \Gamma (M,{\W} ) 
\stackrel{\delta} \to {\mathcal A}^1(M,{\W} ) 
\stackrel{\delta} \to {\mathcal A}^2(M,{\W} ) \stackrel{\delta} \to \ldots, 
 \end{equation}
where 
${\mathcal A}^n (M, {\W})$ denotes  $C^\infty$--sections 
of $n$--form bundle with values in the bundle $\W$,
\[
{\cal A}^k (M,{\W})= {\Gamma} 
(M, \Lambda^k {\mathcal T}^*M \otimes {\W}). 
\]
This complex is exact since $\delta$ is homotopic to identity  by $\delta^*$. 
For each $i$ the  equation for $r_i$  has the form 
 \begin{equation}
 \label{eq:r}
\delta (r_i) = \text{function}(\nabla, r_1, \ldots, r_{i-1}).
 \end{equation}
However it is not difficult to show that  this function is in the kernel
of $\delta$ hence $r_i$ exists.

The  noncommutative associative structure on the Weyl bundle
determines a $\ast$--product
on functions by  the correspondence \ref{eq:ast}.

In fact, the equation $D^2 = 0$
is just the Maurer--Cartan equation for a flat connection. 
One can see the analogy with 
Kazhdan connection \cite{GKF} on the algebra of formal vector fields.
Notice that  $\delta = dx^k
 \frac{\partial }{\partial e^k}$  is of degree $-1$. The flatness of
the connection is given by the  recurrent  procedure, namely starting from the 
terms of degree $-1$ and $0$ one can get other terms step by step. 
 While Kazhdan connection does not have a parameter involved it has the
same structure -- it starts with known $-1$ and $0$ degree terms. Other
terms are of higher degree and can be recovered one by one.

Let us also
mention  here that the connection $D$ can be written as a sum of two terms
-- one is a derivation along the manifold, the usual differential $d$, and
the other is an endomorphism of a fibre, let us call it $\Gamma $. 
Since all endomorphisms
are inner one can write it as an adjoint action with respect to the Moyal 
product. $\Gamma $ acts 
adjointly by an 
operator from $\Gamma (M,{\cal W})$ to $\Gamma (M, {\mathcal T}^*M
\otimes {\cal W})$. 
 \begin{equation}
\label{eq:con}
D = d + \Gamma = d + \dispt[{\gamma}, \cdot]_\circ,
 \end{equation}
where the Lie bracket is the commutator
bracket, and  ${\gamma} \in \Gamma (M, {\mathcal T}^*M \otimes {\cal W})$.
 Then the  equation $D^2 = 0$ becomes:
 \[
d \Gamma + \frac{1}{2} [\Gamma,\Gamma ]_\circ =0.
\]
The same equation for  ${\gamma}$ then is as follows:
 \begin{equation} 
\label{eq:eq}
\omega + d \gamma + \dispt \frac{[\gamma,\gamma]_\circ}{2} = 0,
 \end{equation}
where $\omega$ is a central 2--form.
 This equation states that $D^2$ is given by an
adjoint action of a central  element, so it is zero.
 However it turns out to be very important which exactly form  $\omega$ is 
given in the center by the connection $D$.
 Inner automorphisms of the Weyl algebra are given by the 
adjoint action by  elements of the algebra (Lemma \ref{l:inner}).
Its central extension gives the whole algebra. 
 Curvature of Fedosov connection is zero, however its lift to the
central extension is nonzero.
\begin{definition}
\label{def:char}
The characteristic class of the deformation is the cohomology 
class of the form  $[\omega] \in H^2  (M)[[t]].$
\end{definition}

\section{Calculations from diagrams.}

Fedosov quantization produces the system of recursively defined
equations in order to find  the flat connection and then
another system for the flat sections 
of this connection. 
Fedosov proved that there are no obstructions to solutions.
 In what follows we  show how to obtain these
equations step by step from diagrams, which present all equations at
once. Since there are other situations when one has to solve systems of
recursive equations on graded objects 
we also hope that our presentation might be useful
in some other calculations possibly of completely different origin.

\subsection{Fedosov connection.\\}
\label{chap:A1}

Let ${\cal A} ^p (W^i) =
 \Gamma (M,{\Lambda}^p {\mathcal T}^*M \otimes {\cal W}^i)$.
Let us represent the action of 
 \[
D = - \delta + \nabla + r_1 + r_2 + r_3 + r_4 + \ldots
\]  
 by arrows pointing in directions corresponding to the degree of 
each component.

Namely, $r_k: {\cal A} ^0 (W^i) \to {\cal A} ^1 (W^{i+k})$ is drawn to
go from the point corresponding to the level $i$ in the first column to
the point in the second column $k$ rows down: $\delta $ goes up one row,
$\nabla$ is on the same level, $r_1$ goes down one level and so on.
Same operators act between first and second column.

\vspace{3mm}

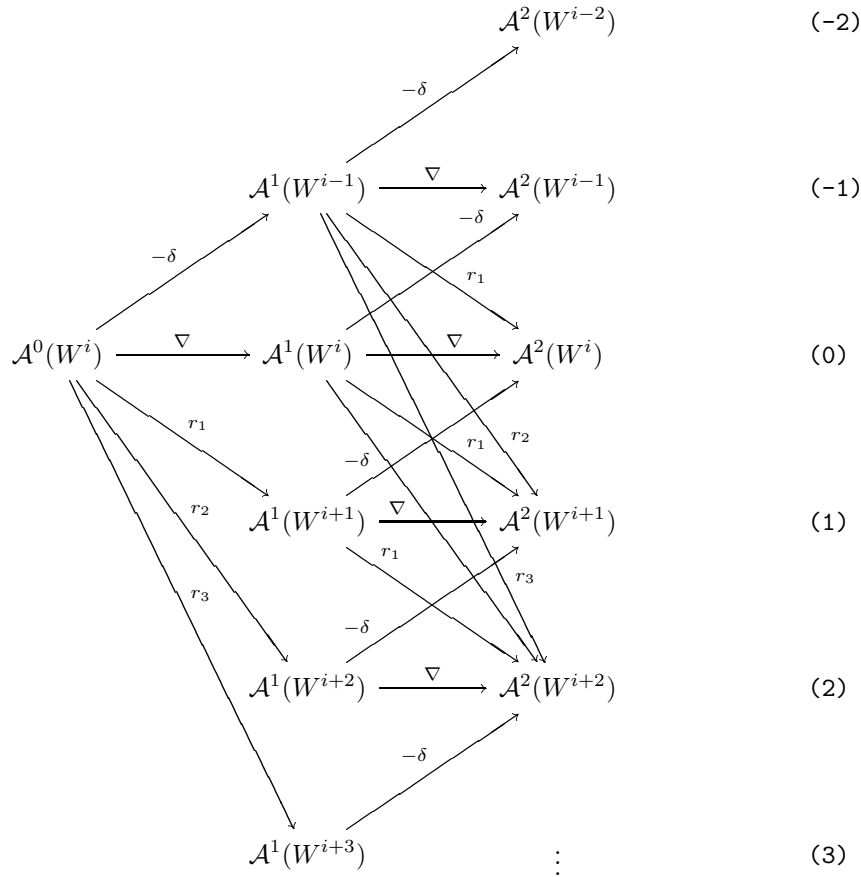
\divide\dgARROWLENGTH by2
\begin{flushleft}
%
% Catcode hack to get typewriter `\' inside arg of another command
% where \verb is illegal.
\begingroup \catcode`|=0 \catcode`\\=12
   |gdef|bbb{{|tt}}%
|endgroup
%
% tighten it up a bit to fit on page in 12pt
\setlength{\dgARROWLENGTH}{1.5em}%
\noindent
$\begin{diagram}
\dgARROWLENGTH=4em
\dgARROWPARTS=6  
  \node[3]{{\cal A}^2(W^{i-2})} 
				\node{\makebox[0pt][l]{\tt\bbb (-2) }}\\
    \node[2]{{\cal A}^1 (W^{i-1})} 
	 \arrow{ne,t}{-\delta}
	\arrow{e,t}{\nabla}
	  \arrow{se,t,4}{r_1}
    	   \arrow{sse,t,5}{r_2} 
	    \arrow{ssse,t,5}{r_3}
    \node{{\cal A}^2(W^{i-1})} 
				\node{\makebox[0pt][l]{\tt\bbb (-1) }}\\
	\node{{\cal A}^0 (W^{i})} 
	    \arrow{ne,t}{-\delta}
	     \arrow{e,t}{\nabla} 
	      \arrow{se,t}{r_1}
		\arrow{sse,t}{r_2}
		 \arrow{ssse,t}{r_3}
 	\node{{\cal A}^1 (W^{i})} 
	  \arrow{ne,t,5}{-\delta}
 	   \arrow{e,t,4}{\nabla} 
	    \arrow{se,t,4}{r_1} 
	      \arrow{sse,t,5}{}
	\node{{\cal A}^2 (W^{i})}
				\node{\makebox[0pt][l]{\tt\bbb (0)  }}\\
	    \node[2]{{\cal A}^1 (W^{i+1})} 
	      \arrow{ne,t,1}{-\delta}
		\arrow{e,t,1}{\nabla}
		  \arrow{se,t,1}{r_1}
	    \node{{\cal A}^2(W^{i+1})} 
				\node{\makebox[0pt][l]{\tt\bbb (1) }}\\    
		\node[2]{{\cal A}^1 (W^{i+2})} 
		   \arrow{ne,t,1}{-\delta}
			\arrow{e,t}{\nabla}
		\node{{\cal A}^2(W^{i+2})}
				\node{\makebox[0pt][l]{\tt\bbb (2) }}\\
		      \node[2]{{\cal A}^1 (W^{i+3})} 
		          \arrow{ne,t}{-\delta} \arrow{e,!}
		      \node{\vdots}
				\node{\makebox[0pt][l]{\tt\bbb (3) }}\\ 
\end{diagram}$
\end{flushleft}

\vspace{-10mm}

{\bf Key observation.}
 The  curvature is equal to $0$ if 
the connection applied  twice to any element 
 $ a \in \Gamma (M,{\cal W})$: $D^2 a = 0$ is $0$ in each degree. In other words:
the sum of arrows coming to the second column 
${\cal A} ^2 (W^i)$ should be  $0$ for every $i$.

Showing that this is true for any  element in 
 ${\cal A} ^0 (W^i) =\Gamma (M,{\cal W}^i)$ for any $i$ will do.

First two terms, $\delta$ of degree $-1$  and $\nabla $ 
of degree $0$, are known,  
our purpose is to find the other terms recursively.

\noindent For every degree $i$  we get  equations on operators:

\begin{itemize}
\item Level (-2):
$\delta^2 = 0$
\item Level (-1):
$- [\delta,\nabla] = 0$
\item Level (0):
$-[\delta,r_1] + {\nabla}^2 = 0$
\item  Level (1):
$-[\delta,r_2] + [\nabla,r_1] = 0$
\item  Level (2):
$-[\delta,r_3] + [\nabla,r_2] + \frac{[r_1,r_1]}{2} = 0$
\item  Level (3):
$-[\delta,r_4] + [\nabla,r_3] + [r_1,r_2] = 0$
\end{itemize}

\noindent  These equations are the graded components of 
$$ D^2 = ( - \delta +
\nabla + [r, .] ) ^2 = \nabla ^2 - [\delta, r] + [\nabla, r] + r^2 = 0$$ 
 where $r =
r_1 + r_2 + \ldots $ It is solved recursively: in each degree $k \geq 0
$ one gets an equation involving only $r_i$ with $ i \leq k$. 

\noindent  Let us show what happens in the first few equations.

{\bf Degree $-2$}. The equation is $\delta ^2 = 0$. It is satisfied 
by  the Lemma (\ref{l:delta}).

{\bf Degree $-1$}. Next one is  $[\delta,\nabla] = 0$. 
It is true by a simple calculation.

\noindent  For this equation we need that the connection $\nabla $ is torsion--free.

{\bf Degree $0$}.  Here is the first nontrivial  calculation.
We have to find such $r_1$ that
$ - [\delta,r_1] = {\nabla}^2 $.

\noindent a) Existence. First of all:
$$[\delta , \nabla ^2] =
[\delta, \nabla] \nabla - \nabla [\delta, \nabla] $$
which is $0$ by the previous equation.

\noindent There is an operator $\delta^*$ which is a homotopy for $\delta$.
\begin{equation}
	\delta^* \delta ^* = 0 , \ \ 
	\delta \delta^* + \delta^* \delta = id \ c
					\label{eq:delta}
\end{equation}
This $c$ is a number of $y$ 's and $dx$ 's, for example for
a term $y^{i_1}\ldots y^{i_p} dx ^ {j_1} \ldots dx ^ {j_q}$
this number $c = p+q $. Let us put $$r_1 =  \delta ^* \nabla ^2 $$
then indeed:
$$ [\delta, r_1] = \delta (r_1) =\delta (\delta ^* \nabla ^2) = \nabla^2$$
and also $\delta ^* r_1 = 0$.

\noindent b) Uniqueness.

\noindent Let $r'_1 = r_1 + \alpha $, such that $\delta^* \alpha = 0$. Then 
$\alpha = \delta ^* \beta$ for some $\beta$. Hence,  $\delta \delta ^* \beta 
= 0$, because $ \delta (r_1 + \alpha ) = \delta r_1$.

\noindent From (\ref{eq:delta})  we get that $\beta = \delta^* \delta \beta$
and $\alpha = \delta^* \beta = \delta ^* (\delta^* \delta \beta) = 0$

\vspace{1cm}

{\bf Degree $1$}.  $[\delta, r_2]  = [\nabla, r_1] $ gives the equation on 
the operator $r_2$.

\noindent a) Existence. Again we show
$$[\delta, [ \nabla, r_1]] = [ [\delta, \nabla], r_1] - 
[ \nabla, [ \delta, r_1]] = 0$$

\noindent b) Uniqueness. $r_2 = \delta ^* ([\nabla, r_1])$ similar to the previous one.

\noindent Recursively getting similar  equations for $r_n$
one finds the Fedosov connection. Here are first few terms:
\[
D = \delta + \nabla + \delta^{-1} \nabla^2 +
 \delta^{-1} \{\nabla, \delta^{-1} \nabla^2 \} + \ldots
\]

%%%%%%%%%%%%%%%%%%%%%%%%%%%%%%%%%%%%%%%%%%

\subsection{Flat sections and the $\ast$--product.\\}

Series in $t$ with 
functional coefficients  are in 
one-to-one correspondence with flat sections of Fedosov connection:
 \[
Q: {C}^\infty(M)[[t]] \to {{\cal A}^0(W)}
 \]
Fedosov connection maps in a unique way each  series
 \[
{a} = a_0 + t a_1 + t^2 a_2+ \ldots
 \]
 to a flat section 
of the Weyl algebra bundle 
 \[
A = {a} + A_1 + A_2 + A_3 + \ldots
 \]
verifying an equation:
 \[
\delta^{-1} (A - { a}) = 0.
\] 
 This last condition makes the operator $Q^{-1}:{{\cal A}^0(W)} \to 
 {C}^\infty(M)[[t]]$  simple, namely it is just an evaluation of
$A \in {{\cal A}^0(W)}$ at zero value of coordinates along the fibres $W$.
This condition could be changed for any other condition fixing the 
zero section in ${{\cal A}^0(W)}$ (see \cite{EW}).

The condition of flatness:
 \[ 
DA = 0
\]
can be represented by the fact that for all $i$ sum of   operators $r_k$
which get to 
${\cal A}^1(W^{i})$ must be $0$.
It again gives a recursive system of equations.

\vspace{3mm}

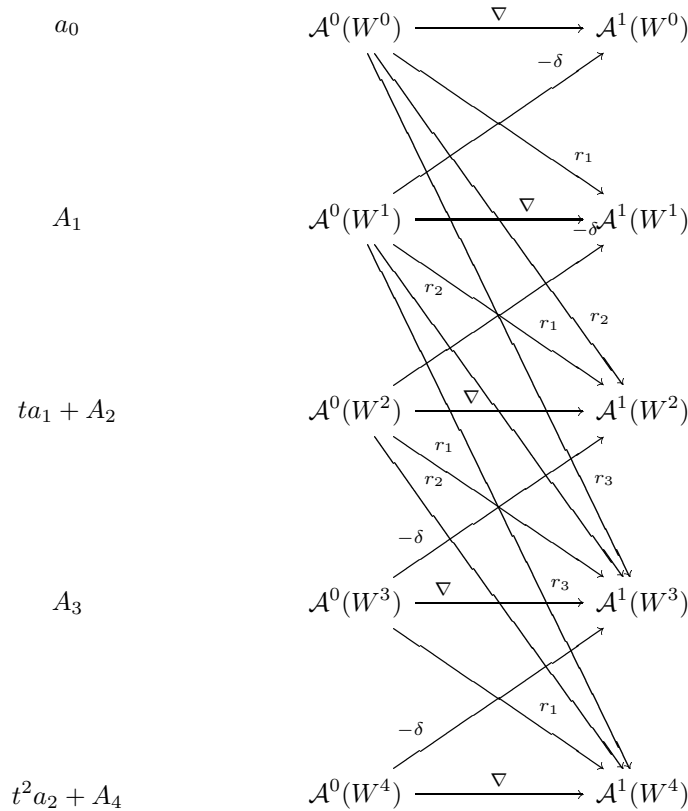
\divide\dgARROWLENGTH by2

\setlength{\dgARROWLENGTH}{1.5em}%
\noindent
$\begin{diagram}
\dgARROWLENGTH=6em
\dgARROWPARTS=6
\node{a_0} 
\arrow{e,!}
	\node{{\cal A}^0(W^{0})} 
	\arrow{e,t}{\nabla} 
	      \arrow{se,t,5}{r_1}
		\arrow{sse,t,5}{r_2}
		 \arrow{ssse,t,5}{r_3}
 	\node{{\cal A}^1 (W^{0})}\\  
\node{A_1} 
 \arrow{e,!}
	\node{{\cal A}^0 (W^{1})} 
	  \arrow{ne,t,5}{-\delta}
 	   \arrow{e,t,4}{\nabla} 
	    \arrow{se,t,4}{r_1} 
	      \arrow{sse,t,1}{r_2}
		\arrow{ssse,t,4}{r_3}
	\node{{\cal A}^1 (W^{1})}\\
\node{ t a_1 +A_2} 
\arrow{e,!}
  	\node{{\cal A}^0 (W^{2})} 
	      \arrow{ne,t,6}{-\delta}
		\arrow{e,t,2}{\nabla}
		  \arrow{se,t,1}{r_1}
			\arrow{sse,t,1}{r_2}
	    \node{{\cal A}^1(W^{2})}\\ 
\node{A_3} 
	   \arrow{e,!}		
\node{{\cal A}^0 (W^{3})} 
		   \arrow{ne,t,1}{-\delta}
			\arrow{e,t,1}{\nabla}
			\arrow{se,t,4}{r_1}
		\node{{\cal A}^1(W^{3})}\\
\node{ t^2 a_2 + A_4} 
	   \arrow{e,!}
		      \node{{\cal A}^0 (W^{4})} 
		          \arrow{ne,t,1}{-\delta} 
				\arrow{e,t} {\nabla}	
		      \node{{\cal A}^1 (W^4)}
\end{diagram}$

\vspace{5mm}

We notice 
that all $r_k$ kill functions, because $r_k$ acts as adjoint operators and
functions are in the center of ${\cal A}^0 (W)$, so $r_k (a_i) = 0 $.
Hence first few  equations following from the diagram above are 

\begin{enumerate}
\item $ \nabla a_0 - \delta A_1 = 0$ 

\item  $ \nabla A_1 - \delta (t a_1 + A_2)= 0$ 

\item $  r_1 A_1 + \nabla (t a_1 +A_2) - \delta A_3   = 0$

\item $ r_2 A_1 + r_1 A_2+ \nabla A_3
				-\delta (t^2 a_2 + A_4) = 0$

\item $ r_3 A_1 + r_2 A_2 + r_1 A_3
				+\nabla (t^2 a_2 + A_4) - \delta A_5= 0$ 
\end{enumerate}

Let $dx^k$ be a local frame in ${\mathcal T}^*M$. Then let the 
 corresponding generators in $W$ 
be $\{e^k\}$. Then $A_i$ are of the 
form $A_{k_1 \ldots k_i}e^{k_1} \ldots e^{k_i}$.

The symplectic connection $\nabla$ locally can be written as:
\[
\nabla  = dx^l \frac{\partial}{\partial x^l} + 
\frac{i}{2 t} [\Gamma_{jkl}dx^le^k e^l, \ \ ]
\]
So for the first terms of the flat section corresponding to $ {   a } = a_0$
we get:
\begin{enumerate}
\item 
$A_1 = \delta^{-1} \nabla  {   a } 
= \delta ^{-1}d{  a} = \partial_l  {   a } e^l$
\item
$A_2 = \delta^{-1} \nabla A_1 
= \delta^{-1} \nabla (\partial_l  {   a } e^l)
= \{\partial_k \partial_l  {   a } 
+ {\Gamma_{kl}}^l\partial_j  {   a } \}e^k e^l$

\end{enumerate}

The  first few terms in the $\ast$--product of two functions 
${   a}, {   b} \in C^\infty(M)$ are:
\[
{   a}*{   b} = {   a }{  b }- \frac{i t}{2}
 \omega^{kl} \partial_i {   a}  \partial_j {  b }
- t^2  (\partial_i \partial_j {  a }+ 
\Gamma_{kl}^l \partial_l {  a }) \omega^{im} \omega^{jn}
(\partial_m \partial_n  {  b }+ 
\Gamma_{mn}^k \partial_k {  b }) + \ldots 
\]

%%%%%%%%%%%%%%%%%%%%%%%%%%%%%%%%%%%%%%%%%
\subsection{Standard example of deformation quantization.\\}

The procedure of deformation quantization requires calculations which
are not obvious and most of the time do not give nice formulas. However
in few cases for particular manifolds one can calculate  the
$\ast$--product  explicitly. The first trivial example 
is the quantization of ${\R}^{2n}$ with a standard symplectic form. 
Let $\{ x^1, \ldots, x^{2n} \}$ be 
a local coordinate system at some point $x \in {\R}^{2n}$.
The Darboux symplectic form in these coordinates is
\[
\omega = dx^i \wedge dx^{i+n}, \quad 1 \leq i \leq n
\]

The
standard symplectic form and the  trivial connection $\nabla =d$
 give an algebra of
differential operators in ${\R}^{2n}$. 
Namely, using calculations from the diagrams 
we get the  Fedosov connection to be:
\[ 
D = d - \delta \quad \text{or in coordinates} \quad D = dx^k 
(\frac{\partial}{\partial x^k} - \frac{\partial}{\partial e^k}).
\]

Flat section of such a connection corresponding to a 
function $ {a}$ under the quantization map is as follows
\[
A = { {a}} + e^k \frac{\partial{  a}}{\partial x^k} + 
e^k e^l \frac{{\partial^2}{  a}}{\partial x^k \partial x^l}
+ \cdots.
\]
 We see that it gives a formula for Taylor decomposition  of a function 
$ {a}$ at a point $x$. In fact the $e^k$ terms can be considered 
as jets.  Then the $\ast$--product of two flat sections 
is given by the formula (\ref{eq:Moyal}). It is easy to deduce that 
for two functions $ {a}$ and $ {b}$ the $\ast$--product is 
\begin{equation}
\begin{split}
{  a}\ast {  b} & =
\exp \  \{- i t \frac{\partial}
{\partial y^k} \frac{\partial  }{\partial z^{n+k}}\}{  a (y)} 
{  b (z)} |_{y=z=x}\\
  &={  a}{  b} -
{i t }\frac{\partial{  a}}{\partial x^{k}}
\frac{\partial{  b}}{\partial x^{n+k}} + \frac {t^2 }{2}
(\frac{\partial^2{  a}}{\partial x^{k}\partial x^{l}})
(\frac{\partial^2{  b}}{{\partial x^{n+k}}{\partial x^{n+l}}}) + \cdots.
\end{split}
\label{eq:R2n}
\end{equation}
Let us   map $C^\infty ({\R}^{2n})$
 to differential operators 
on ${\R}^{n}$, considered as polynomials on $T^*{\R}^{n}$. Then the 
 $\ast$--product gives exactly the product of differential 
symbols.

\begin{remark}
The same scheme actually works for any cotangent  
bundle ${\mathcal T}^*M$ with the 
canonical symplectic form -- the
quantized algebra of functions on  ${\mathcal T}^*M$ is 
isomorphic to the algebra of
differential operators on $M$. This observation leads to various types of 
index theorems \cite{Gutt}, also \cite{Ezra},\cite{NT}.
\end{remark}
%%%%%%%%%%%%%%%%

\end{document}